\documentclass[12pt]{article}
\usepackage{amsmath,latexsym,epsf}
\usepackage{amsfonts}
\usepackage{amssymb}
\usepackage{mathrsfs} 
\usepackage{amsmath}
\usepackage[english]{babel} 
\usepackage[latin1]{inputenc}
\usepackage[all]{xy}
\begin{document}

\newcommand{\nc}{\newcommand}
\def\PP#1#2#3{{\mathrm{Pres}}^{#1}_{#2}{#3}\setcounter{equation}{0}}
\def\ns{$n$-star}\setcounter{equation}{0}
\def\nt{$n$-tilting}\setcounter{equation}{0}
\def\Ht#1#2#3{{{\mathrm{Hom}}_{#1}({#2},{#3})}\setcounter{equation}{0}}
\def\qp#1{{${(#1)}$-quasi-projective}\setcounter{equation}{0}}
\def\mr#1{{{\mathrm{#1}}}\setcounter{equation}{0}}
\def\mh#1{{{\mathcal{#1}}}\setcounter{equation}{0}}
\newcommand{\LL}{\ell\ell\,}
\newcommand{\edge}{\ar@{-}}
\newcommand{\bxy}{\xymatrix}
\newtheorem{Th}{Theorem}[section]
\newtheorem{Def}[Th]{Definition}
\newtheorem{Lem}[Th]{Lemma}
\newtheorem{Pro}[Th]{Proposition}
\newtheorem{Cor}[Th]{Corollary}
\newtheorem{Rem}[Th]{Remark}
\newtheorem{Exm}[Th]{Example}
\def\Pf#1{{\noindent\bf Proof}.\setcounter{equation}{0}}
\def\>#1{{ $\Rightarrow$ }\setcounter{equation}{0}}
\def\<>#1{{ $\Leftrightarrow$ }\setcounter{equation}{0}}
\def\bskip#1{{ \vskip 20pt }\setcounter{equation}{0}}
\def\sskip#1{{ \vskip 5pt }\setcounter{equation}{0}}
\def\midskip#1{{ \vskip 10pt }\setcounter{equation}{0}}
\def\bg#1{\begin{#1}\setcounter{equation}{0}}
\def\ed#1{\end{#1}\setcounter{equation}{0}}
\def\r#1{{\rm{#1}}}


\title{\bf Finitistic Dimension Conjecture and Conditions on Ideals}
\smallskip
\author{{\small Jiaqun WEI} {\thanks {Supported by the National Science Foundation of China
(No.10601024)}}\\
\small Department of Mathematics,
Nanjing Normal University \\
\small Nanjing 210097, P.R.China\\ \small Email:
weijiaqun@njnu.edu.cn}
\date{}
\maketitle
\baselineskip 16.3pt
%
%
\begin{abstract}
\vskip 10pt%
The notion of Igusa-Todorov classes is introduced in connection
with the finitistic dimension conjecture. As application we
consider conditions on special ideals which imply the
Igusa-Todorov and other finiteness conditions on modules proving
the finitistic dimension conjecture and related conjectures in
those cases.
\midskip\
%
%

\end{abstract}
%
\vskip 30pt
%

{\noindent \Large\bf Introduction}\vskip 15pt
 \hskip 15pt

Throughout the paper, we work on artin algebras and finitely
generated left modules.

Let $A$ be an artin algebra. Recall that the  finitistic dimension
of $A$, denoted by $\mr{findim}A$, is defined to be the supremum
of the projective dimensions of all finitely generated modules of
finite projective dimension. The famous   finitistic dimension
conjecture claims that $\mr{findim}A<\infty$, for any artin
algebra $A$. The conjecture is still open and is related to many
other homological conjectures (e.g., the Gorenstein symmetry
conjecture, the Wakamatsu-tilting conjecture and the generalized
Nakayama conjecture).

Until recently, only a few classes of algebras were known to have
finite finitistic dimension, see for instance [\ref{GKK},
\ref{GZh}, \ref{Hp}] etc.. In particular, the finitistic dimension
conjecture holds for algebras of representation dimension at most
3, algebras with radical cube zero, monomial algebras and left
serial algebras, etc..

 All the above-mentioned algebras are Igusa-Todorov algebras,
a notion introduced in [\ref{Wit}] for algebras whose finiteness
of finitistic dimensions is assured by the Igusa-Todorov functor
[\ref{ITf}]. Many algebras which are recently proved to have
finite finitistic dimension are also Igusa-Todorov [\ref{Wit}].

In this paper, we introduce the notion of Igusa-Todorov classes
and apply them to obtain new Igusa-Todorov algebras. The main
results are as follows.

\vskip  15pt

\bg{Th}\label{Mth}
Let $A$ be an artin algebra.

$(1)$ Assume that $I,J$ are two ideals of $A$ such that
$IJ\mr{rad}A=0$. Then $A$ is Igusa-Todorov provided one of the
following conditions is satisfied.

\hskip 30pt $(i)$ $A/I$ is representation-finite,
$\mr{pd}_AJ<\infty$ and $A/J$ is syzygy-finite.

\hskip 30pt $(ii)$ $A/J$ is representation-finite,
$\mr{pd}_AI<\infty$ and $A/I$ is syzygy-finite.

\hskip 30pt $(iii)$ Both $A/I$ and $A/J$ are syzygy-finite and
both $_AI$ and $_AJ$ have finite projective dimension.

$(2)$ Assume that $0=I_0\subseteq I_1\subseteq\cdots\subseteq I_n$
are ideals in $A$ such that $\mr{pd}_AI_{i+1}<\infty$ and
$(I_{i+1}/I_i)\mr{rad}(A/I_i)=0$, for each $0\le i\le n-1$. If
$A/I_n$ is syzygy-finite (resp., syzygy-bounded, Igusa-Todorov),
then $A$ is also  syzygy-finite (resp., syzygy-bounded,
Igusa-Todorov).

In particular, the finitistic dimension of $A$ is finite in all
the above cases.
\ed{Th}%

%

Recall now the following well-known conjectures (see for instance
[\ref{Hp}, \ref{MR}]). Assume that $A$ is an artin algebra.

\vskip  15pt
{\bf Gorenstein symmetry conjecture}\ \  $\mr{id}(_AA)<\infty$ if
and only if $\mr{id}(A_A)<\infty$, where id denotes the injective
dimension.

{\bf Wakamatsu-tilting conjecture}\ \  Let $_A\omega$ be a
Wakamatsu-tilting module. (1) If $\mr{pd}_A\omega<\infty$, then
$\omega$ is tilting. (2) If $\mr{id}_A\omega<\infty$, then
$\omega$ is cotilting.

{\bf Generalized Nakayama conjecture}\ \  Each indecomposable
injective $A$-module occurs as a direct summand in the minimal
injective resolution of $_AA$.

\vskip   15pt

Note that the Gorenstein symmetry conjecture and the generalized
Nakayama conjecture are special cases of the second
Wakamatsu-tilting conjecture. Moreover, if the finitistic
dimension conjecture holds for $A$ and $A^o$ (hereafter $A^o$
denotes the opposite algebra of $A$), then all the above
conjectures hold.
%
%

Small's  results [\ref{Sm}] shows that if $A$ has a nilpotent
ideal $I$ such that $\mr{pd}_AI<\infty$ and $(A/I)^o$ has finite
finitistic dimension, then $A^o$ has also finite finitistic
dimension. Hence, combining with the fact that, for a
syzygy-finite algebra $A$, both the finitistic dimensions of $A$
and $A^o$ are finite (see for instance [\ref{Wab}]), we have the
following corollary.

\vskip  15pt

\bg{Cor}\label{Mc}%
The Gorenstein symmetry conjecture, the Wakamatsu-tilting
conjecture and the generalized Nakayama conjecture hold for
algebras in {\rm Theorem \ref{Mth}} $(1)$, provided that moreover,
one of there considered ideals is nilpotent.
\ed{Cor}%

%

In particular, Small's results [\ref{Sm}] shows that the
finitistic dimension of $A^o$ is finite provided that
$\mr{pd}_A\mr{rad}^2A<\infty$ or $\mr{pd}_A\mr{rad}^3A<\infty$
(see for instance [\ref{Zho}]). It is natural to ask if the
finitistic dimension of $A$ is also finite in these cases. The
following corollary gives a partial answer to that question.

%

\bg{Cor}\label{Mcr}
Let $A$ be an artin algebra. Then the finitistic dimension of $A$
is finite provided that $A$ satisfies one of the following two
conditions.

$(1)$ $\mr{rad}^5A=0$ and $\mr{pd}_A\mr{rad}^2A<\infty$.

$(2)$ $\mr{rad}^nA=0$ and $\mr{pd}_A\mr{rad}^iA<\infty$ for each
$i\ge 3$.

In particular, the Gorenstein symmetry conjecture, the
Wakamatsu-tilting conjecture and the generalized Nakayama
conjecture hold for algebras in two  cases.
\ed{Cor}

%

In the last section, we give examples to show that our results do
apply to obtain algebras with finite finitistic dimension. In
particular, our methods are compared with the others, for
instance, the recent results of Huards, Lanzilotta and Mendoza
[\ref{HLM}] using their infinite radical.

%
\vskip 30pt
\section{Preliminaries}

\hskip 15pt

Throughout this paper, $A$ is always an artin algebra. We denote
by $A\mr{-mod}$ the category of all (finitely generated)
$A$-modules. If $M$ is an $A$-module, we use $\mr{pd}_AM$ to
denote the projective dimension of $_AM$ and use $\Omega^i_AM$ to
denote the $i$-th syzygy of $M$ (we set $\Omega^0_AM:=M$). The
notion $\mr{add}_AM$ denotes the category consists of all direct
summands of finite direct sums of $M$.

Let $\mh{C}$ be a class of $A$-modules. We denote
$\mr{gd}\mh{C}:=sup\{\mr{pd}_AM|M\in\mh{C}\}$. We use
$\Omega^n\mh{C}$ to denote the class of $A$-modules which are
isomorphic to $\Omega^n_AM$ for some $M\in\mh{C}$. The class
$\mh{C}$ is said to be {\it add-finite} provided that
$\mh{C}\subseteq \mr{add}_AC$ for some $A$-module $C$, or
equivalently,  the number of non-isomorphic indecomposable direct
summands of objects in $\mh{C}$ is finite. Moreover, we say that
$\mh{C}$ is ($n$-){\it syzygy-finite} provided that there is some
$n$ such that $\Omega^n\mh{C}$ is add-finite. If $A\mr{-mod}$ is
syzygy-finite, then we also say that $A$ is syzygy-finite. Note
that if $\mh{C}$ is $n$-syzygy-finite, then $\Omega^i\mh{C}$ is
add-finite for any $i\ge n$ and $\Omega^i\mh{C}$ is
$(n-i)$-syzygy-finite for each $1\le i< n$.

More generally, we say that a class of $A$-modules $\mh{C}$ is
{\it add-bounded} if the non-isomorphic indecomposable direct
summands of objects in $\mh{C}$ have bounded length. We say that
$\mh{C}$ is ($n$-){\it syzygy-bounded} if there exists $n$ so that
$\Omega^iC$ is add-bounded for each $i\ge n$. If $A\mr{-mod}$ is
syzygy-bounded, then we also say that $A$ is syzygy-bounded.
%
Obviously an add-finite  (resp., a syzygy-finite) class is
add-bounded (resp., syzygy-bounded).

Note that  syzygy-bounded classes have finitistic projective
dimension by a well-known theorem of Jensen and Lenzing [\ref{JL}]
which is based on a theorem of Schofield [\ref{Scf}].

 It is easy to see that if $\mh{D}\subseteq\mh{C}$ and
$\mh{C}$ is add-finite (resp.,  syzygy-finite, add-bounded,
syzygy-bounded), then $\mh{D}$ is also  add-finite (resp.,
syzygy-finite,
add-bounded, syzygy-bounded). 

%
%

Now let us recall the definition of Igusa-Todorov algebras
introduced in [\ref{Wit}]. By definition  an artin algebra $A$ is
called ($n$-){\it Igusa-Todorov} provided that there exists a
fixed $A$-module $_AV$ and a nonnegative integer $n$ such that,
for any $M\in A\mr{-mod}$, there is an exact sequence $0\to V_1\to
V_0\to \Omega^n_AM\to 0$ with $V_0, V_1\in \mr{add}_AV$. A
remarkable property of Igusa-Todorov algebras, which is also the
reason why Igusa-Todorov algebras were named, is that these
algebras can be shown to satisfy the finitistic dimension
conjecture by using the following results on the Igusa-Todorov
functor defined in [\ref{ITf}].

\bg{Lem}\label{fclem}
{\rm [\ref{ITf}]} For any artin algebra $A$, there is a functor
$\Psi$ which is defined on the objects of $A\mr{-mod}$ and takes
non-negative integers as values, such that

$(1)$ $\Psi(M)=\mr{pd}_AM$ provided that $\mr{pd}_AM<\infty$.

$(2)$ $\Psi(X)\le \Psi(Y)$ whenever $\mr{add}_AX\subseteq
\mr{add}_AY$. The equation holds in case $\mr{add}_AX=
\mr{add}_AY$.

$(3)$ If $0\to X\to Y\to Z\to 0$ is an exact sequence in
$A\mr{-mod}$, then $\mr{pd}_AZ\le \Psi(X\oplus Y)+1$, provided
that $\mr{pd}_AZ<\infty$.
\ed{Lem}%

The class of Igusa-Todorov algebras is large and includes all
algebras of representation dimension at most 3, algebras with
radical cube zero, monomial algebras and syzygy-finite algebras,
etc. [\ref{Wit}].

%
\vskip 30pt
\section{Igusa-Todorov and other finiteness conditions}


\hskip 15pt


Inspired by the definition of Igusa-Todorov algebras, we introduce
here the following notions.

\bg{Def}\label{ITdef}%
Let $\mh{C}$ be a class of $A$-modules. Then $\mh{C}$ is said to
be $(n$-$)$Igusa-Todorov provided that there is an $A$-module
$_AV$ and some nonnegative integer $n$ such that, for any $M\in
\Omega^n\mh{C}$, there is an exact sequence $0\to V_1\to V_0\to
M\to 0$ with $V_0, V_1\in \mr{add}_AV$. The module $_AV$ is then
called an $(n$-$)\mh{C}$-Igusa-Todorov module.
\ed{Def}%

By definition, we immediately obtain that all syzygy-finite
classes are Igusa-Todorov.

\bg{Rem}\label{rem}%
 $(1)$ If $\mh{C}\subseteq\mh{D}$ and $\mh{D}$ is
Igusa-Todorov, then $\mh{C}$ is Igusa-Todorov too.

$(2)$ $A$ is an Igusa-Todorov algebra if and only if $A\mr{-mod}$
is Igusa-Todorov.

$(3)$ $\mh{C}$ is $n$-Igusa-Todorov if and only if
$\Omega^n\mh{C}$ is $0$-Igusa-Todorov.

$(4)$ If $\Omega^i\mh{C}$ is Igusa-Todorov for some $i\ge 1$, then
$\mh{C}$ is also Igusa-Todorov.

$(5)$ In general, there are many $n$-$\mh{C}$-Igusa-Todorov
modules for an $n$-Igusa-Todorov class $\mh{C}$.
\ed{Rem}


\bg{Lem}\label{oITl}
Let $\mh{C}$ be a class of $A$-modules. If $\mh{C}$ is
$0$-Igusa-Todorov, then $\Omega^i\mh{C}$ is also
$0$-Igusa-Todorov, for any $i\ge 1$.
\ed{Lem}%

\Pf. Let $_AV$ be an $n$-$\mh{C}$-Igusa-Todorov module. Then for
any $M\in \mh{C}$, there is an exact sequence $0\to V_1\to V_0\to
M\to 0$ with $V_0, V_1\in \mr{add}_AV$. By the Horseshoe Lemma,
for any $\Omega^i_AM$, there is an exact sequence
$0\to\Omega^i{_AV_1}\to\Omega^i{_AV_0}\oplus P\to\Omega^i_AM\to 0$
for some projective $A$-module $P$. Hence $\Omega^i\mh{C}$ is
$0$-Igusa-Todorov with a $0$-$\Omega^i\mh{C}$-Igusa-Todorov module
$\Omega^i{_AV}\oplus A$. \hfill $\Box$

%
%
%
%


Immediately by Lemma \ref{oITl} and Remark \ref{rem} (3), we
obtain the following proposition.

\bg{Pro}\label{ITnm}%
Let $\mh{C}$ be a class of $A$-modules. If $\mh{C}$ is
$n$-Igusa-Todorov, then $\Omega^i\mh{C}$ is also $n$-Igusa-Todorov
for any $i\ge 1$. In particular  $\mh{C}$ is $m$-Igusa-Todorov,
for any $m\ge n$ in that case.
\ed{Pro}

\vskip 15pt

The following result gives a method to obtain Igusa-Todorov
classes.

\bg{Pro}\label{Mlem}
Let $\mh{C}, \mh{D}, \mh{E}$ be three classes of $A$-modules such
that, for any $E\in\mh{E}$, there is an exact sequence $0\to C\to
D\to E\to 0$ with $C\in\mh{C}$ and $D\in\mh{D}$.

$(1)$ If $\mh{C}, \mh{D}$ are syzygy-finite, then $\mh{E}$ is
Igusa-Todorov.

$(2)$ If $\mh{C}$ is syzygy-finite (resp., syzygy-bounded) and
$\mr{gd}\mh{D}<\infty$, then $\mh{E}$ is also syzygy-finite
(resp., syzygy-bounded).

$(3)$ If $\mh{C}$ is Igusa-Todorov and $\mr{gd}\mh{D}<\infty$,
then $\mh{E}$ is also Igusa-Todorov.
\ed{Pro}%

\Pf. By assumption, for any $E\in\mh{E}$, there is an exact
sequence $0\to C\to D\to E\to 0$ with $C\in\mh{C}$ and
$D\in\mh{D}$. Thus, for any $n\ge 1$, we obtain an exact sequence
$0\to\Omega^nC\to\Omega^nD\oplus P\to\Omega^nE\to 0$ for some
projective $A$-module $P$, by the Horseshoe Lemma.

(1) Take $n$ to be an integer such that both $\Omega^n\mh{C}$ and
$\Omega^n\mh{D}$ are add-finite, and denote that $U=\oplus U_i$
(resp., $V=\oplus V_i$), where $U_i$ (resp., $V_i$) lists all the
non-isomorphic indecomposable direct summands of objects in
$\Omega^n\mh{C}$ (resp., $\Omega^n\mh{D}$). Then we have that
$\Omega^nC\in\mr{add}_AU$ and $\Omega^nD\in\mr{add}_AV$. It
follows that $\mh{E}$ is $n$-Igusa-Todorov with an
$n$-$\mh{E}$-Igusa-Todorov module $U\oplus V\oplus A$.

(2) Take $n$ to be an integer such that $\Omega^n\mh{C}$ is
add-finite (resp., add-bounded) and $\mr{gd}\mh{D}\le n$. Then
$\Omega^nD$ is projective. It follows that
$\Omega^nC\simeq\Omega^{n+1}E\oplus Q$ for some projective $Q$.
Now we easily obtain that $\Omega^{n+1}\mh{E}$ is add-finite
(resp., add-bounded). Hence, $\mh{E}$ is syzygy-finite (resp.,
syzygy-bounded) by the definition involved.

(3) Take $n$ to be an integer such that $\mh{C}$ is
$n$-Igusa-Todorov and $\mr{gd}\mh{D}\le n$. Similarly to the proof
of (2), we obtain that $\Omega^nC\simeq\Omega^{n+1}E\oplus Q$ for
some projective $Q$. Note that there is an exact sequence $0\to
V_1\to V_0\to \Omega^nC\to 0$ with $V_0, V_1\in \mr{add}_AV$,
where $_AV$ is an $n$-$\mh{C}$-Igusa-Todorov module. Since $Q$ is
projective, we can also obtain an exact sequence $0\to V_1'\to
V_0'\to \Omega^{n+1}E\to 0$ with $V_0', V_1'\in \mr{add}_AV$. It
follows that $\mh{E}$ is $(n+1)$-Igusa-Todorov with $_AV$ an
$(n+1)$-$\mh{E}$-Igusa-Todorov module.
 \hfill $\Box$

\vskip 15pt

\bg{Cor}\label{Mcor}%
Let $\mh{C}, \mh{D}, \mh{E}$ be three classes of $A$-modules such
that, for any $E\in\mh{E}$, there is an exact sequence $0\to C\to
D_0\to\cdots\to D_n\to E\to 0$ with $C\in\mh{C}$ and each
$D_i\in\mh{D}$. If $\mh{C}$ is syzygy-finite (resp.,
syzygy-bounded, Igusa-Todorov) and $\mr{gd}\mh{D}<\infty$, then
$\mh{E}$ is also syzygy-finite (resp., syzygy-bounded,
Igusa-Todorov).
\ed{Cor}%

\Pf.  Denote $\mh{E}_0:=\mh{C}_0$, and inductively,
$\mh{E}_{i+1}:=\{M|$ there is an exact sequence $0\to C\to D\to
M\to 0$ with $C\in\mh{E}_i$ and $D\in\mh{D}\}$. Then by
assumptions and Proposition \ref{Mlem}, inductively we obtain that
each $\mh{E}_i$ is syzygy-finite (resp., syzygy-bounded,
Igusa-Todorov). Note that $\mh{E}\subseteq\mh{E}_{n+1}$, so
$\mh{E}$ is also syzygy-finite (resp., syzygy-bounded,
Igusa-Todorov).
%
 \hfill $\Box$

%

%
\vskip 15pt
%
%
%
%

We say that two classes of $A$-modules $\mh{C,D}$ are {\it syzygy
equivalent} if there exist $n,m$ so that
$\mr{add}(\Omega^n\mh{C})=\mr{add}(\Omega^m\mh{D})$ (up to the
class $\mr{add}_AA$). For example, if there is a finite exact
sequence in which all but two modules $M,N$ have finite projective
dimension then $\mr{add}M, \mr{add}N$ are syzygy equivalent. It is
obvious that for a class $\mh{C}$, all classes $\Omega^i\mh{C}$
and $\mh{C}$ are syzygy-equivalent. It is also easy to see that
$\mh{C}$ is syzygy-finite (resp., syzygy-bounded, Igusa-Todorov)
if and only it is in the syzygy equivalent class of an add-finite
(resp., add-bounded, $0$-Igusa-Todorov) class.

We have the following result.

\vskip 15pt

\bg{Pro}\label{syeqhom}%
Let $A\to R$ be a homomorphism of of artin algebras. If
$\mr{pd}_AR<\infty$, then the restriction functor: $$ R-\mr{mod}
\rightarrow A-\mr{mod} $$ preserves syzygy equivalence and sends
syzygy equivalence classes of modules to syzygy equivalence
classes of modules.
\ed{Pro}

\Pf. For a class $\mh{C}$ of $R$-modules, we denote by
$\mh{_{A}{C}}$ the image of $\mh{C}$ under the restrict functor.
Note that for any $R$-module $M\in\Omega\mh{C}$ (resp.,
$N\in\mh{C})$, there is an exact sequence $0\to M\to P\to N\to 0$
with $P\in\mr{add}_RR$ and $N\in\mh{C}$ (resp.,
$M\in\Omega\mh{C}$). Applying the restrict functor (which is
clearly exact), we obtain that, for any $A$-module
${_A}M\in\mh{_A}\Omega\mh{C}$ (resp., ${_A}N\in\mh{_AC})$, there
is an exact sequence of $A$-modules $0\to {_A}M\to {_A}P\to
{_A}N\to 0$ with $P\in\mr{add}_AR$ and ${_A}N\in\mh{C_A}$ (resp.,
${_A}M\in\mh{_A}\Omega\mh{C}$). Since $\mr{pd}_AR<\infty$, we see
that $\mh{_A}(\Omega\mh{C}), \mh{_AC}$ are syzygy equivalent. It
follows that $\mh{_A}(\Omega^i\mh{C}), \mh{_AC}$ are also syzygy
equivalent for any $i\ge 0$.

Assume now that $\mh{C,D}\subseteq R-\mr{mod}$ are syzygy
equivalent. It is enough to show that $\mh{{_A}C,{_A}D}\subseteq
A-\mr{mod}$ are syzygy equivalent.

By the definition, there exist $n,m$ so that
$\mr{add}(\Omega^n\mh{C})=\mr{add}(\Omega^m\mh{D})$ up to the
class $\mr{add}_RR$. Restricting them to $A\mr{-mod}$, we have
that
$\mr{add}\mh{_A}(\Omega^n\mh{C})=\mr{add}\mh{_A}(\Omega^m\mh{D})$
up to the class $\mr{add}_AR$. Since $\mr{pd}_AR<\infty$, it is
easy to see that $\mh{_A}(\Omega^n\mh{C}),
\mh{_A}(\Omega^m\mh{D})$ are syzygy equivalent. Note that
$\mh{_A}(\Omega^n\mh{C})$ and $\mh{_AC}$ (resp.,
$\mh{_A}(\Omega^m\mh{D})$ and $\mh{_AD}$) are syzygy equivalent by
the above argument, so we obtain that $\mh{{_A}C,{_A}D}$ are
syzygy equivalent.
 \hfill$\Box$

\vskip 15pt

\bg{Cor}\label{Mcor2}%
Let $A\to R$ be a homomorphism of artin algebras such that
$\mr{pd}_AR<\infty$.
%
If the class of $R$-modules $\mh{C}$ is syzygy-finite (resp.,
syzygy-bounded, Igusa-Todorov), then its image under the restrict
functor, $\mh{_{A}{C}}$, is also syzygy-finite (resp.,
syzygy-bounded, Igusa-Todorov). In addition, if  $\mh{D}\subseteq
A\mr{-mod}$ is syzygy equivalent to a subclass of $\mh{_AC}$, then
$\mh{D}$ is also syzygy-finite (resp., syzygy-bounded,
Igusa-Todorov).
%
%
%
\ed{Cor}%

\Pf. Note that the restrict functor $R\mr{-mod}\to A\mr{-mod}$ is
exact and length preserving, so all relevant finiteness conditions
(e.g., add-finite, add-bounded, $0$-Igusa-Todorov) on a class of
$R$-modules are inherited by its image in $A\mr{-mod}$. Hence, if
$\mh{C}$ is syzygy-finite (resp., syzygy-bounded, Igusa-Todorov),
then $\mh{_{A}{C}}$ is also syzygy-finite (resp., syzygy-bounded,
Igusa-Todorov) by Proposition \ref{syeqhom}. The remained part is
obvious.
 \hfill$\Box$
%

\vskip 15pt

Let $A$ and $R$ be artin algebras. Following [\ref{xi1}], we say
that $R$ is a left idealized extension of $A$, provided that
$A\subseteq R$ has the same identity and $\mr{rad}A$ is a left
ideal in $R$. It was shown that every finite dimensional algebra
admits a chain of left idealized extensions which is of finite
length and ends with a representation-finite algebra [\ref{xi1},
Lemma 4.6].

Generalizing [\ref{xi2}, Theorem 3.1], it is proved in [\ref{Wit},
Theorem 2.9] that if $A=A_0\subseteq\cdots\subseteq A_m=R$ are
artin algebras such that, for each $0\le i\le m-1$, $A_{i+1}$ is a
left idealized extension of $A_i$ and $\mr{pd}_AA_{i+1} <\infty$,
and $R$ is Igusa-Todorov, then the finitistic dimension of $A$ is
finite. One part of the following result shows that the artin
algebra $A$ obtained in that way is in fact Igusa-Todorov.

\bg{Th}\label{rid}%
Assume that $A=A_0\subseteq\cdots\subseteq A_m=R$ are artin
algebras such that, for each $0\le i\le m-1$, $A_{i+1}$ is a left
idealized extension of $A_i$ and $\mr{pd}_AA_{i+1} <\infty$. If
$R$ is  syzygy-finite (resp., syzygy-bounded, Igusa-Todorov), then
$A$ is also syzygy-finite (resp., syzygy-bounded, Igusa-Todorov).
In particular, the finitistic dimension of $A$ is finite.
%
%
\ed{Th}

\Pf. Take any $X_0\in A\mr{-mod}$. By [\ref{xi1}, Erratum, Lemma
0.1], if $A_1$ is a left idealized extension of $A$, then for any
$X\in A\mr{-mod}$, $\Omega_A^2(X)\simeq \Omega_{A_1}Z\oplus Q$ as
$A_1$-modules for some $Z,Q\in A_1\mr{-mod}$ with $_{A_1}Q$
projective. Thus, by assumptions and inductively, we see that, for
each $0\le i\le m-1$ and each $X_i\in A_i\mr{-mod}$,
$\Omega_{A_i}^2(X_i)\simeq \Omega_{A_{i+1}}(X_{i+1})\oplus
P_{i+1}$ as $A_{i+1}$-modules for some $X_{i+1}, P_{i+1}\in
A_{i+1}\mr{-mod}$ and $_{A_{i+1}}P_{i+1}$ projective. Hence we
have the following exact sequence of $A_{i+1}$-modules, where
$P'_{i+1}$ is the projective cover of $\Omega_{A_{i+1}}(X_{i+1})$:
$$0\to \Omega^2_{A_{i+1}}(X_{i+1})\to P'_{i+1}\oplus P_{i+1}\to
\Omega_{A_i}^2(X_i)\to 0.$$

It is easy to see that all these exact sequences can be restricted
to the exact sequences of $A$-modules. Thus we can combine all
these exact sequences and obtain the following exact sequence of
$A$-modules:
$$0\to \Omega^2_{A_{m-1}}(X_{m-1})\to P'_{m-1}\oplus P_{m-1}\to\cdots\to P_1'\oplus P_1\to
\Omega_A^2(X_0)\to 0.$$

 Denote
now $\mh{D}=\mr{add}_A(\oplus_{i=0}^{m} A_i)$, then each
$P_i'\oplus P_i\in \mh{D}$ and $\mr{gd}\mh{D}<\infty$ by
assumption. Let $\mh{R}$ be the class of all $A$-modules which are
restricted from $R$-modules.  Then we obtain that $A\mr{-mod}$ is
syzygy equivalent to a subclass  of $\mh{R}$ from the above exact
sequence. Therefore, if that $R$ is syzygy-finite (resp.,
syzygy-bounded, Igusa-Todorov), then $A\mr{-mod}$ is also
syzygy-finite (resp., syzygy-bounded, Igusa-Todorov), by Corollary
\ref{Mcor2}.
%
 \hfill $\Box$

%
\vskip 15pt

\bg{Cor}\label{ridc}%
Assume that $A=A_0\subseteq\cdots\subseteq A_m=R$ are artin
algebras such that, for each $0\le i\le m-1$, $A_{i+1}$ is a left
idealized extension of $A_i$ and $\mr{pd}_{A_i}A_{i+1} <\infty$.
If $R$ is syzygy-finite (resp., syzygy-bounded, Igusa-Todorov) and
$R^o$ has finite finitistic dimension, then the Gorenstein
symmetry conjecture, the Wakamatsu-tilting conjecture and the
generalized Nakayama conjecture hold for $A$.
\ed{Cor}

\Pf. Since $\mr{pd}_{A_i}A_{i+1} <\infty$ for each $0\le i\le
m-1$, we easily obtain that $\mr{pd}_{A}A_{i+1} <\infty$. It
follows from Theorem \ref{rid} that the finitistic dimension of
$A$ is finite. Moreover, by assumptions and [\ref{Wfd}, Corollary
3.7], the finitistic dimension of $A^o$ is also finite. Now the
conclusion follows from these two facts.
 \hfill $\Box$

%
\vskip 30pt

\section{Conditions on ideals}

\hskip 15pt

The following result can be compared with [\ref{xi1}, Theorem
3.2], where it was proved that if $A$ has two ideals $I,J$ such
that $IJ=0$ and $A/I$ and $A/J$ are representation-finite, then
the finitistic dimension of $R$ is finite. In fact, it is shown in
[\ref{Wit}] that $A$ is $0$-Igusa-Todorov in that case.

\bg{Th}\label{ij}%
Assume that $I,J$ are two ideals of $A$ such that $IJ\mr{rad}A=0$.
Then $A$ is Igusa-Todorov provided one of the following conditions
is satisfied.

$(1)$ $A/I$ is representation-finite, $\mr{pd}_AJ<\infty$ and
$A/J$ is syzygy-finite.

$(2)$ $A/J$ is representation-finite, $\mr{pd}_AI<\infty$ and
$A/I$ is syzygy-finite.

$(3)$ Both $A/I$ and $A/J$ are syzygy-finite and both $_AI$ and
$_AJ$ have finite projective dimension.

In particular, the finitistic dimension of $A$ is finite in all
three cases.
\ed{Th}%

\Pf. Denote by $\mh{I_A}$ (resp., $\mh{J_A}$) the class of all
$A$-modules which are also $A/I$-modules (resp.,  $A/J$-modules).
For any $N\in \Omega(A\mr{-mod})$, we have that $N\subseteq
\mr{rad}_AP$ for some projective $_AP$. Denote $C_1:=JN$ and
$C_2:=N/JN$, then we have an exact sequence $0\to C_1\to N\to
C_2\to 0$. Note that, by assumptions, $IC_1=IJN\subseteq
IJ\mr{rad}_AP=0$ and $JC_2=J(N/JN)=0$ , so $C_1$ is also an
$A/I$-module and $C_2$ is also an $A/J$-module. Now consider the
following pullback diagram, where the right column is taken such
that $_AP$ is the projective cover of $_AC_2$.

\bskip\

 \setlength{\unitlength}{0.09in}
 \begin{picture}(50,22)

 \put(27,1){\makebox(0,0)[c]{$0$}}
 \put(35,1){\makebox(0,0)[c]{$0$}}

                 \put(27,4){\vector(0,-1){2}}
                 \put(35,4){\vector(0,-1){2}}

 \put(12,5){\makebox(0,0)[c]{$0$}}
                             \put(14,5){\vector(1,0){2}}
 \put(18,5){\makebox(0,0)[c]{$C_1$}}
                             \put(20,5){\vector(1,0){2}}
 \put(27,5){\makebox(0,0)[c]{$N$}}
                             \put(31,5){\vector(1,0){2}}
 \put(35,5){\makebox(0,0)[c]{$C_2$}}
                             \put(37,5){\vector(1,0){2}}
 \put(41,5){\makebox(0,0)[c]{$0$}}

                 \put(18,9){\line(0,-1){2}}
                 \put(17.5,9){\line(0,-1){2}}
                 \put(27,9){\vector(0,-1){2}}
                 \put(35,9){\vector(0,-1){2}}

 \put(12,11){\makebox(0,0)[c]{$0$}}
                             \put(14,11){\vector(1,0){2}}
 \put(18,11){\makebox(0,0)[c]{$C_1$}}
                             \put(20,11){\vector(1,0){2}}
 \put(27,11){\makebox(0,0)[c]{$C_1\oplus P$}}
                             \put(31,11){\vector(1,0){2}}
 \put(35,11){\makebox(0,0)[c]{$P$}}
                             \put(37,11){\vector(1,0){2}}
 \put(41,11){\makebox(0,0)[c]{$0$}}

                 \put(27,14){\vector(0,-1){2}}
                 \put(35,14){\vector(0,-1){2}}

 \put(27,16){\makebox(0,0)[c]{$\Omega_AC_2$}}
                             \put(30,16.5){\line(1,0){2}}
                             \put(30,16){\line(1,0){2}}
 \put(35,16){\makebox(0,0)[c]{$\Omega_AC_2$}}

                 \put(27,19){\vector(0,-1){2}}
                 \put(35,19){\vector(0,-1){2}}

 \put(27,21){\makebox(0,0)[c]{$0$}}
 \put(35,21){\makebox(0,0)[c]{$0$}}

\end{picture}
\bskip\

Note that $\Omega_AC_2\in\Omega{\mh{J_A}}$ and $C_1\oplus
P\in\mh{D}:=\{C \oplus P| C \in\mh{I_A}, P\in\mr{add}_AA\}$, so if
we show that both $\Omega{\mh{J_A}}$ and $\mh{D}$ are
syzygy-finite, then $A$ is Igusa-Todorov by Remark \ref{rem} (3),
Proposition \ref{Mlem} and the middle column in the diagram.

$(1)$ If $\mr{pd}_AJ<\infty$ and $A/J$ is syzygy-finite, then
$\mr{pd}_AA/J<\infty$. By Corollary \ref{Mcor2}, we obtain that
the class $\mh{J_A}$ is syzygy-finite. It follows that
$\Omega{\mh{J_A}}$ is also syzygy-finite. Since $A/I$ is
representation-finite, the class $\mh{D}$ is clearly also
syzygy-finite.

$(2)$ Similarly as in (1), we have that $\mh{I_A}$ is
syzygy-finite. It follows that $D$ is also syzygy-finite. Since
$A/J$ is representation-finite, the class $\Omega{\mh{J_A}}$ is
clearly also syzygy-finite.

$(3)$ Similarly as in (1), we have that both $\mh{I_A}$ and
$\mh{J_A}$ are syzygy-finite by Corollary \ref{Mcor2}.
Consequently, both $\Omega{\mh{J_A}}$ and $\mh{D}$ are
syzygy-finite.
 \hfill $\Box$

%
\vskip 15pt

The following list mentions some known classes of syzygy-finite
algebras.

$\bullet$ Algebras of finite representation type (in fact
$0$-syzygy-finite).

$\bullet$ Algebras of finite global dimension.

$\bullet$ Monomial algebras (indeed $2$-syzygy-finite)
[\ref{Zhi}].

$\bullet$ Left serial algebras (indeed $2$-syzygy-finite)
[\ref{Zhs}].

$\bullet$ Torsionless-finite algebras  (indeed $1$-syzygy-finite)
c.f. [\ref{R}], including:

\hskip 30pt $\circ$ Algebras $A$ with $\mr{rad}^nA = 0$ and
$A/\mr{rad}^{n-1}A$  representation-finite.

\hskip 30pt $\circ$ Algebras $A$ with $\mr{rad}^2A = 0$.

\hskip 30pt $\circ$ Minimal representation-infinite algebras.

\hskip 30pt $\circ$ Algebras stably equivalent to hereditary
algebras.

\hskip 30pt $\circ$ Right glued algebras and left glued algebras.

\hskip 30pt $\circ$ Algebras of the form $A/\mr{soc}A$ with $A$ a
local algebra of quaternion type.

\hskip 30pt $\circ$  Special biserial algebras.

$\bullet$ Algebras possessing a left idealized extension which is
torsionless-finite (indeed $2$-syzygy-finite) [\ref{Wit}].
%

%
\vskip 15pt

\bg{Cor}\label{i2r}%
Let $A$ be an artin algebra.  If there is an ideal $I$ such that
$I^2\mr{rad}A=0\ ($for instance $I^2=0)$, $\mr{pd}_AI<\infty$ and
$A/I$ is syzygy-finite, then the finitistic dimension of $A$ is
finite.
\ed{Cor}

%

A special case of the above corollary is as follows.

\bg{Cor}\label{radn}%
Let $A$ be an artin algebra with $\mr{rad}^{2n+1}A=0$. If
$\mr{pd}_A\mr{rad}^nA<\infty$ and $A/\mr{rad}^nA$ is
syzygy-finite, then the finitistic dimension of $A$ is finite.
\ed{Cor}

%

In particular, we obtain the following result.

\bg{Cor}\label{rad2}%
Let $A$ be an artin algebra with $\mr{rad}^5A=0$. If
$\mr{pd}_A\mr{rad}^2A<\infty$, then the finitistic dimension of
$A$ is finite.
\ed{Cor}

\Pf. Note that $A/\mr{rad}^2A$ is an algebra with radical square
zero and hence is torsionless-finite. Thus the conclusion follows
from Corollary \ref{radn}.
 \hfill $\Box$

\vskip 15pt

\bg{Pro}\label{ir}%
Let $A$ be an artin algebra and $I$ be an ideal with
$I\mr{rad}A=0$ and $\mr{pd}_AI<\infty$. If $A/I$ is  syzygy-finite
(resp., syzygy-bounded, Igusa-Todorov), then  $A$ is also
syzygy-finite (resp., syzygy-bounded, Igusa-Todorov). In
particular, the finitistic dimension of $A$ is finite.
%
%
\ed{Pro}%

\Pf. By assumptions, for any $N\in A\mr{-mod}$, we have  that
$I\Omega_AN\subseteq I\mr{rad}(P_N)=0$, where $P_N$ is the
projective cover of $N$. Hence $\Omega_AN\in\mh{I_A}$, where
$\mh{I_A}$ denotes the class of all $A$-modules which are also
$A/I$-modules.  It follows that $A\mr{-mod}$ is syzygy equivalent
to a subclass  of $\mh{I_A}$. Moreover, we have that
$\mr{pd}_AA/I<\infty$, since $\mr{pd}_AI<\infty$ by the
assumption. Hence, if $A/I$ is syzygy-finite (resp.,
syzygy-bounded, Igusa-Todorov), then $A\mr{-mod}$ is also
syzygy-finite (resp., syzygy-bounded, Igusa-Todorov) by Corollary
\ref{Mcor2}.
%
%
 \hfill $\Box$

\vskip 15pt

Generally, we have the following result.

\bg{Th}\label{in}%
Assume that $0=I_0\subseteq I_1\subseteq\cdots\subseteq I_n$ are
ideals in $A$ such that $\mr{pd}_AI_{i+1}<\infty$ and
$(I_{i+1}/I_i)\mr{rad}(A/I_i)=0$, for each $0\le i\le n-1$. If
$A/I_n$ is  syzygy-finite (resp., syzygy-bounded, Igusa-Todorov),
then $A$ is also syzygy-finite (resp., syzygy-bounded,
Igusa-Todorov). In particular, the finitistic dimension of $A$ is
finite.
%
\ed{Th}

\Pf. Denote that $A_i=A/I_i$ for  each $0\le i\le n$. So $A_0=A$.
Note that $\mr{pd}_AA_i<\infty$ for each  $0\le i\le n$, since
$\mr{pd}_AI_i<\infty$ by assumption.

Take any $X_0\in A_0\mr{-mod}$ and denote $X_1=\Omega_{A_0}X_0$,
then we have an exact sequence $0\to X_1\to P_{X_0}\to X_0\to 0$,
where $P_{X_0}$ is the projective cover of $X_0$. Since
$I_1X_1=I_1\Omega_{A_0}X_0\subseteq I_1\mr{rad}_A(P_{X_0})=0$ by
assumptions, we see that $X_1$ is also an $A_1$-module.
Inductively, for each $0\le i\le n-1$ and each $A_i$-module $X_i$,
we obtain an exact sequence $0\to X_{i+1}\to P_{X_i}\to X_i\to 0$,
where $P_{X_i}$ is the projective cover of the $A_i$-module $X_i$,
such that $X_{i+1}$ is also an $A_{i+1}$-module.

Now restricting these exact sequences to $A$-modules and combining
them, we obtain an exact sequence $0\to X_n\to P_{X_{n-1}}\to
\cdots\to P_{X_0}\to X_0\to 0$ of $A$-modules, where $X_n$ is also
an $A_n$-module and each $P_{X_i}$ is projective as an
$A_i$-module.

Denote now $\mh{D}=\mr{add}_A(\oplus_{i=0}^{n-1} A_i)$, then each
$P_{X_i}\in \mh{D}$ and $\mr{gd}\mh{D}<\infty$ by assumption. Let
$\mh{A}_n$ be the class of all $A$-modules which are also
$A_n$-modules. Then we obtain that $A\mr{-mod}$ is syzygy
equivalent to a subclass  of $\mh{A}_n$ from the above exact
sequence. Now the conclusion easily follows from Corollary
\ref{Mcor2}.
 \hfill $\Box$

\vskip 15pt

\bg{Cor}\label{rn}%
Let $n$ be the integer such that $\mr{rad}^nA=0$.

$(1)$ If $\mr{pd}_A\mr{rad}^iA<\infty$ for each $i\ge 3$, then $A$
is Igusa-Todorov.

$(2)$ If $\mr{pd}_A\mr{rad}^iA<\infty$ for each $i\ge 2$, then $A$
is syzygy-finite.

In particular, the finitistic dimension of $A$ is finite in either
case.
\ed{Cor}

\Pf. Consider the ideals $0=\mr{rad}^nA\subseteq
\mr{rad}^{n-1}A\subseteq\cdots\subseteq\mr{rad}^3A\subseteq\mr{rad}^2A$.
Then clearly
$(\mr{rad}^{i-1}A/\mr{rad}^iA)\cdot\mr{rad}(A/\mr{rad}^iA)=0$.

(1)  Since $A/\mr{rad}^3A$ is an algebra with radical cube zero,
it is Igusa-Todorov by [\ref{Wit}]. Hence we obtain the conclusion
by assumptions and Theorem \ref{in}.

 (2) Since $A/\mr{rad}^2A$ is an algebra with radical square zero,
it is syzygy-finite. Again we obtain the conclusion by assumptions
and Theorem \ref{in}.
 \hfill $\Box$

\vskip 15pt

As special cases of Corollaries \ref{rad2} and \ref{rn}, we obtain
the following result.

\bg{Cor}\label{r4}%
If $\mr{rad}^4A=0$, then the finitistic dimension of $A$ is finite
provided that $\mr{pd}_A\mr{rad}^2A<\infty$ or
$\mr{pd}_A\mr{rad}^3A<\infty$.
\ed{Cor}


\vskip 30pt

\section{Examples}

\hskip 15pt

In this section, we will present examples to show that our results
do apply to check the finiteness of the finitistic dimensions of
some algebras.

\vskip 15pt

\begin{Exm}\label{exmm} {\rm [\ref{HLM}]} %
Consider the bound quiver algebra $A=kQ/J$ where $Q$ is given by

$$
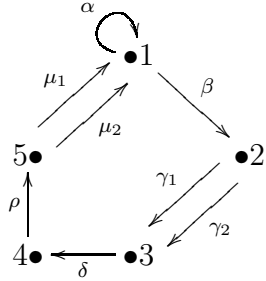
\xymatrix{ & {\bullet}1\ar@(l,u)[]^{\alpha} \ar[dr]^{\beta} & \\
{5}{\bullet}\ar@<1ex>[ru]^{\mu_1} \ar@<-1ex>[ru]_{\mu_2}& &{\bullet}2 \ar@<1ex>[dl]^{\gamma_2} \ar@<-1ex>[dl]_{\gamma_1}\\
{4}{\bullet}\ar[u]^{\rho} & {\bullet}{3} \ar[l]^{\delta} & }$$

\noindent and $J$ is generated by the set of paths
$\{\alpha^3 , \beta\alpha, \alpha\mu_i\rho,\beta\mu_i
,\delta\gamma_1-\delta\gamma_2\}$
where $i=1,2$. Then $A$ is syzygy-finite. In particular, the
finitistic dimension of $A$ is finite.
\end{Exm}

The algebra was proved to have finite finitistic dimension in
[\ref{HLM}] by showing that $\LL^\infty A=3$. Here the notion
$\LL^\infty A$ denotes the infinite-layer length of $A$, see
[\ref{HLM}] for details.  In fact, one main result in [\ref{HLM}]
asserts that an artin algebra $A$ has finite finitistic dimension
whenever $\LL^\infty A\le 3$, generalizing the known result for
algebras with radical cube zero.

Now we give a different proof of the finiteness of the finitistic
dimension of $A$. Moreover, we prove that $A$ is syzygy-finite.

To see this, let $I_1$ be the ideal generated by $\beta$ and $I_2$
be the ideal generated by $\beta, \gamma_1,\gamma_2, e_2$. Then
 $I_1\subseteq I_2$ and
$I_1\mr{rad}A=0=(I_2/I_1)\mr{rad}(A/I_1)$. Moreover,
$\mr{pd}_AI_1<\infty$ and $\mr{pd}_AI_2<\infty$. Note that $A/I_2$
is of the form

$$
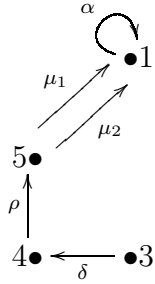
\xymatrix{ & {\bullet}1\ar@(l,u)[]^{\alpha} & \\
{5}{\bullet}\ar@<1ex>[ru]^{\mu_1} \ar@<-1ex>[ru]_{\mu_2}& &\\
{4}{\bullet}\ar[u]^{\rho} & {\bullet}{3} \ar[l]^{\delta} & }$$

\noindent with relations
$0=\alpha^3=\alpha\mu_i\rho$,
where $i=1,2$. Clearly $A/I_2$ is a monomial algebra and hence is
syzygy-finite. Now by Theorem \ref{in} we obtain that $A$ is also
syzygy-finite.

%
%
%

\vskip 15pt

\begin{Exm}\label{exmm2}
Consider the bound quiver algebra $A=kQ/J$ where $Q$ is the same
as in Example \ref{exmm}, but $J$ is generated by the set of paths
$\{\alpha^6 , \beta\alpha, \alpha\mu_i\rho,\beta\mu_i
,\delta\gamma_1-\delta\gamma_2\}$,
where $i=1,2$. Then $A$ is syzygy-finite. In particular, the
finitistic dimension of $A$ is finite.
\end{Exm}

The algebra is revisited from Example \ref{exmm} and the argument
is the same. Note that here $\LL^\infty A=6$, so we can't use the
results in [\ref{HLM}].

\vskip 15pt

\begin{Exm}\label{xi4e} {\rm [\ref{xi4}]} %
Consider the bound quiver algebra $A=kQ/J$ where $Q$ is given by

$$
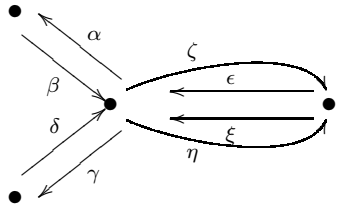
\xymatrix{{\bullet} \ar@<-1ex>[dr]_{\beta} & & & &\\
 & {\bullet} \ \ \ar@(l,u)[urdr]^{\zeta} \ \ \ar@(l,d)[rdru]_{\eta} \ar@<-1ex>[ul]_{\alpha} \ar@<1ex>[dl]^{\gamma} &  &{\bullet} \ar@<1ex>[ll]^{\xi} {\ar@<-1ex>[ll]_{\epsilon}} \\
{\bullet} \ar@<1ex>[ru]^{\delta} & & & }$$

\noindent and $J$ is generated by the set of paths
$\{\alpha\beta, \alpha\epsilon, \alpha\delta, \alpha\xi,
\gamma\beta, \gamma\delta, \gamma\epsilon, \gamma\xi,
\zeta\epsilon, \epsilon\zeta, \xi\zeta, \eta\epsilon, \eta\xi,
\beta\alpha-\delta\gamma\}$.
Then $A$ is syzygy-finite. In particular, the finitistic dimension
of $A$ is finite.
\end{Exm}

It was shown in [\ref{xi4}] that $A$ has a semisimple extension
$A'$ such that $\mr{rad}A=\mr{rad}A'$ and $A'$ is a trivially
twisted extension of two algebras of finite finitistic dimension.
Consequently, by results in [\ref{xi4}], the algebra $A$ has
finite finitistic dimension, see [\ref{xi4}] for details.

Now we give a different proof of the finiteness of the finitistic
dimension of $A$. Moreover, we also show that $A$ is
syzygy-finite.

To see this, let $I$ be the ideal generated by $\alpha$. Then
$I\mr{rad}A=0$. Moreover, $\mr{pd}_AI<\infty$. Note that $A/I$ is
of the form

$$
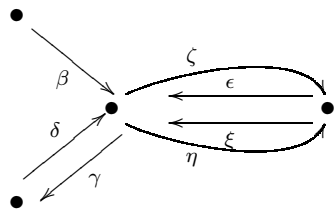
\xymatrix{{\bullet} \ar[dr]_{\beta} & & & &\\
 & {\bullet} \ \ \ar@(l,u)[urdr]^{\zeta} \ \ \ar@(l,d)[rdru]_{\eta} \ar@<1ex>[dl]^{\gamma} &  &{\bullet} \ar@<1ex>[ll]^{\xi} {\ar@<-1ex>[ll]_{\epsilon}} \\
{\bullet} \ar@<1ex>[ru]^{\delta} & & & }$$

\noindent with relations
$0=\delta\gamma=\gamma\beta=
\gamma\delta=\gamma\epsilon=\gamma\xi=\zeta\epsilon=
\epsilon\zeta=\xi\zeta=\eta\epsilon=\eta\xi$.
Clearly $A/I$ is a monomial algebra and hence is syzygy-finite.
Now by Proposition \ref{ir} we obtain that $A$ is also
syzygy-finite.

\vskip 15pt

\begin{Exm}\label{xi4e2}%
Consider the bound quiver algebra $A=kQ/J$ where $Q$ is given by

$$
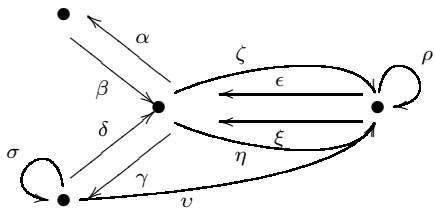
\xymatrix{{\bullet} \ar@<-1ex>[dr]_{\beta} & & & &\\
 & {\bullet} \ \ \ar@(l,u)[urdr]^{\zeta} \ \ \ar@(l,d)[rdru]_{\eta} \ar@<-1ex>[ul]_{\alpha} \ar@<1ex>[dl]^{\gamma} &  &{\bullet} \ar@<1ex>[ll]^{\xi} {\ar@<-1ex>[ll]_{\epsilon}} \ar@(u,r)[]^{\rho} &\\
{\bullet} \ar@<1ex>[ru]^{\delta} \ar@(u,l)[]_{\sigma}
\ar@(l,d)[rrru]_{\upsilon} & & & }$$

\noindent and $J$ is generated by paths
$\{\rho^6, \sigma^6, \alpha\beta, \alpha\epsilon, \alpha\delta,
\alpha\xi, \gamma\beta, \gamma\delta, \gamma\epsilon, \gamma\xi,
\zeta\epsilon, \epsilon\zeta, \xi\zeta, \eta\epsilon, \eta\xi,
\beta\alpha-\delta\gamma\}$.
Then $A$ is syzygy-finite. In particular, the finitistic dimension
of $A$ is finite.
\end{Exm}

The algebra is revisited from Example \ref{xi4e} and the argument
is almost the same.

\vskip 15pt

\begin{Exm}\label{new}%
Consider the bound quiver algebra $A=kQ/J$ where $Q$ is given by

$$
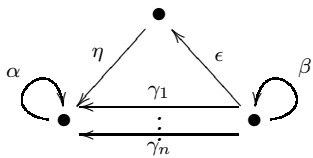
\xymatrix{ & {\bullet} \ar[ld]_{\eta} & \\
{\bullet}  \ar@(l,u)[]^{\alpha} & \vdots & {\bullet}
\ar@(r,u)[]_{\beta}  \ar@<-1ex>[ll]_{\gamma_1}
\ar@<1ex>[ll]^{\gamma_n}   \ar[lu]_{\epsilon}}$$

\noindent and $J$ is generated by the set of paths
$\{\alpha\gamma_i\beta-\alpha\gamma_j\beta$
for all $1\le i<j\le n$, $\alpha^3, \beta^6\}$. Then $A$ is
Igusa-Todorov. In particular, the finitistic dimension of $A$ is
finite.
\end{Exm}

In fact, let $I$ be the ideal generated by $\{\gamma_i, 1\le i\le
n-1\}$. Then $I^2=0$. Moreover, $\mr{pd}_AI<\infty$. Note that
$A/I$ is of the form

$$
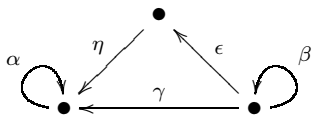
\xymatrix{ & {\bullet} \ar[ld]_{\eta} & \\
{\bullet} \ar@(l,u)[]^{\alpha} & & {\bullet} \ar@(r,u)[]_{\beta}
\ar[ll]_{\gamma}  \ar[lu]_{\epsilon} }$$

\noindent with relations
$0=\alpha\gamma\beta=\alpha^3=\beta^6$.
Clearly $A/I$ is a monomial algebra and hence is syzygy-finite.
Now by Corollary \ref{i2r} we obtain that $A$ is Igusa-Todorov.

Similarly, we can also obtain that $A^o$ is also Igusa-Todorov and
hence has also finite finitistic dimension.

\vskip 15pt

We note that the Gorenstein symmetry conjecture, the
Wakamatsu-tilting conjecture and the generalized Nakayama
conjecture also hold for the algebras in these examples.
\vskip 30pt
\begin{center}
{\bf ACKNOWLEDGEMENTS}
\end{center}
\hskip 15pt The author thanks for the referee's highly valuable
comments and suggestions. In particular, the notions of
``syzygy-bounded" and ``syzygy-equivalent'' are due to the
referee, as well as Propositions \ref{syeqhom}.


{\small

}

%
%
%

\end{document}